%%%%%%%%%%%%%%%%%%%%%%%%%%    parity2.tex    %%%%%%%%%%%%%%%%%%%%%%%%%%
%%%%%%%%%%%%%%%%%%%%%%%%%%    Plain TeX     %%%%%%%%%%%%%%%%%%%%%%%%%%

\input amssym.def
\input amssym.tex 

\font\twbf=cmbx12

%%%%%%       M A C R O S       %%%%%%%

%%%%    CONTROL

\def\bn{\bigskip\noindent}
\def\mn{\medskip\noindent}
\def\sn{\smallskip\noindent}

\def\hb{\hfill\break}
\def\hat{\widehat}
\def\o{\overline}

%%%%%%%     OPERATORS

\def\Aut{{\mathop{\rm Aut}}}
\def\End{{\mathop{\rm End}}}
\def\Hom{{\mathop{\rm Hom}}}

\def\Supp{{{\rm Supp}}}
\def\Tr{{\mathop{\rm Tr}}}
\def\ord{{\mathop{\rm ord}}}
\def\corank{{\mathop{\rm corank}}}
\def\rk{{\mathop{\rm rk}}}
\def\tors{{{\rm tors}}}

\def\s{{S_\infty}}
\def\x{{X_\infty}}

\def\ran #1{{r_{an}(E/#1)}}
\def\sp #1{{s_p(E/#1)}}
\def\Mod #1{{{}_{#1}{\rm Mod}}}

\def\ki{{{K_\infty}}}
\def\gn{{{\Gamma_n}}}
\def\kn{{{K_n}}}

\def\invlim{\mathop{\vtop{\ialign{##\crcr$\hfill{\lim}\hfil$\crcr
\noalign{\kern1pt\nointerlineskip}\leftarrowfill\crcr\noalign
{\kern -3pt}}}}\limits}

\def\dirlim{\mathop{\vtop{\ialign{##\crcr$\hfill{\lim}\hfil$\crcr
\noalign{\kern1pt\nointerlineskip}\rightarrowfill\crcr\noalign
{\kern -3pt}}}}\limits}

\def\iso{\buildrel \sim \over \longrightarrow}
\def\ho{\hookrightarrow}
\def\lo{\longrightarrow}
\def\Lo{\Longrightarrow}
\def\nd{\!\not\kern2.3pt\mid}

\def\({\left(}
\def\){\right)}

%%%%    GENERAL   LETTERS

\def\wp{{{\goth p}}}

\def\Q{{{\bf Q}}}
\def\C{{{\bf C}}}

\def\N{{{\cal N}}}
\def\O{{{\cal O}}}
\def\Z{{{\bf Z}}}
\def\qb{{\overline{\bf Q}}}
\def\qp{{{\bf Q}_p}}
\def\zp{{{\bf Z}_p}}

\font\cyr=wncyr10
\def\sha{\hbox{\cyr X}}

%%%%%%%%%   DIAGRAMS AND ARROWS   %%%%%%%%%

%%%%%%   END  OF    M A C R O S       %%%%%%%

\centerline{\twbf On the parity of ranks of Selmer groups II}

\bn
\centerline{\bf Jan Nekov\'a\v{r}}

\sn
\centerline{\it 11 November 2000}

\vskip20pt
\noindent
Let

$$
E: y^2 = x^3 + Ax + B \qquad\qquad (A,B\in\Q)
$$
be an elliptic curve over $\Q$ of conductor $N$. Thanks to the work of
Wiles and his followers [BCDT] we know that $E$ is modular, i.e. there
exists a non-constant map $\pi : X_0(N) \lo E$ defined over $\Q$ and

$$
L(E,s) = \sum_{n=1}^\infty a_n n^{-s} = L(f,s)
$$
for a normalized newform $f\in S_2(\Gamma_0(N))$.

For a large class of number fields $F$ (which includes all solvable Galois
extensions of $\Q$), this is known to imply that the $L$-function
$L(E/F, s)$ has a holomorphic continuation to $\C$ and a functional
equation relating the values at $s$ and $2-s$. For such $F$, denote by

$$
\ran F := \ord_{s=1} L(E/F,s)
$$
the {\sl analytic rank} of $E$ over $F$.

Over an arbitrary number field $F$, the $m$-descent on $E$ gives rise
(for every integer $m\geq 1$) to the Selmer group $S(E/F,m)$ sitting
in the standard exact sequence

$$
0 \lo E(F) \otimes \Z/m\Z \lo S(E/F,m) \lo \sha(E/F)[m] \lo 0.
$$
Fix a prime number $p$ and put

$$
\eqalign{
S(E/F) &= S_p(E/F) = \dirlim_n S(E/F, p^n)\cr
X(E/F) &= X_p(E/F) = \invlim_n S(E/F, p^n)\cr
\sp F &= \corank_\zp S_p(E/F) = \rk_\zp X_p(E/F).\cr}
$$
The conjecture of Birch and Swinnerton-Dyer predicts that

$$
\ran \Q \buildrel ? \over = \rk_\Z E(\Q).
\leqno (BSD)
$$
As in [NePl], we are interested in a rather weak consequence of (BSD)
(and the conjectural finiteness of the Tate-\v{S}afarevi\v{c} group),
namely the

\mn
{\bf Parity conjecture for Selmer groups:}
$\qquad \ran \Q \buildrel ? \over \equiv \sp \Q \pmod 2 .$

\sn
Our main result is the following

\proclaim Theorem A. Let $E$ be an elliptic curve over $\Q$ with good
ordinary reduction at $p$. Then the parity conjecture
$$
\ran \Q \equiv \sp \Q \pmod 2
$$
holds.

\sn
See [NePl] for a discussion of earlier results in this direction. Our
method of proof is similar to that in [NePl]; the only difference is
that we use anticyclotomic deformations instead of Hida families.
A recently proved conjecture of Mazur ([Maz2], [Co], [Va2]) on non-vanishing
of Heegner points in anticyclotomic $\zp$-extensions plays the role of
Greenberg's conjecture assumed in [NePl].

\vskip24pt
\centerline{\bf 1. Heegner points}

\vskip24pt\noindent
In this section we recall the basic setup of Heegner points on $E$
(see [Gro] for a more detailed account).

\sn
{\bf 1.1} \ Let $K = \Q(\sqrt D)$ be an imaginary quadratic field of
discriminant $D < 0$. We assume that $K$ satisfies the following
``Heegner condition" of Birch [Bi]:

$$
\qquad\qquad{\rm Every\ prime\ } q\vert N {\rm\ splits\ in\ } K.
\leqno (Heeg)
$$
Under this assumption there exists an ideal $\N \subset \O_K$ such that
$\O_K/\N \iso \Z/N\Z$ (of course, $\N$ is not unique; we choose one).

\sn
{\bf 1.2} \ For every integer $c\geq 1$ denote by $\O_c = \Z + c\O_K$
the unique order of $O_K$ of conductor $c$. If $(c,N) = 1$, then 

$$
\N_c := \O_c\cap\N
$$
is an invertible ideal in $\O_c$ satisfying $\O_c/\N_c \iso \Z/N\Z$
(and hence also $\N_c^{-1}/\O_c \iso \Z/N\Z$, since $\N_c$ is
invertible). The cyclic $N$-isogeny

$$
[\C/\O_c \lo \C/\N_c^{-1}]
$$
induced by the identity map on $\C$ defines a non-cuspidal point
on the modular curve $X_0(N)$, which is defined over $H_c$, the ring class
field of conductor $c$ over $K$. The image of this point under the
fixed modular parametrization $\pi : X_0(N) \lo E$ will be denoted
by $\o x_c \in E(H_c)$ -- this is a {\sl Heegner point of conductor
$c$ on $E$}.

\sn
{\bf 1.3} \ From now on, we consider only conductors of the form $c = p^n$
for a fixed prime number $p$. The field $H_{p^\infty} = \bigcup H_{p^n}$
then contains the {\sl anticyclotomic $\zp$-extension} $\ki$ of $K$, which
can be characterized by the following properties: $\ki = \bigcup \kn$,
$G(\kn/K) \iso \Z/p^n\Z$, $\kn/\Q$ is a Galois extension with $G(\kn/\Q)
\iso D_{2p^n}$ (the dihedral group). There is $n_0\geq 0$ such that

$$
H_{p^{n+1}}\cap \ki = K_{n+n_0} \qquad\qquad (n\geq 0).
$$

We use the following standard notation:

$$
\Gamma = G(\ki/K), \qquad \Gamma_n = G(\ki/\kn), \qquad \Lambda =
\zp[[\Gamma]].
$$
Fix an isomorphism $\Gamma\iso\zp$ (i.e. fix a topological generator
$\gamma$ of $\Gamma$).

From now on, we assume that the following condition of ``good reduction"
is satisfied:

$$
p\nd N.
\leqno (G)
$$
This assumption implies that the Heegner points

$$
\o x_{p^{n+1}} \in E(H_{p^{n+1}})
$$
are defined. We put

$$
x_{n+n_0} := \Tr_{H_{p^{n+1}}/K_{n+n_0}} (\o x_{p^{n+1}}) =
\Tr_{H_{p^\infty}/\ki} (\o x_{p^{n+1}}) \in E(K_{n+n_0})
\qquad\qquad (n\geq 0).
$$

\sn
{\bf 1.4} \ In his lecture at ICM 1983, Mazur formulated (among others)
the following

\proclaim Conjecture [Ma2]. \quad $(\exists\, n\geq n_0)
\quad x_n\not\in E(\kn)_\tors$.

\sn
This conjecture was recently proved by Vatsal [Va2] under the assumptions
that $|D|$ is prime and $p$ does not divide the class number of $K$,
and by Cornut [Co] under the assumption (G).
Both [Co] and [Va2] build upon an earlier work [Va1] of Vatsal.

\sn
{\bf 1.5} \ The fundamental distribution relation for Heegner points
([PR], Lemma 2, p.432) states that

$$
\Tr_{K_{n+1}/\kn}(x_{n+1}) = a_p x_n - x_{n-1} \qquad\qquad (n\geq n_0+1).
$$
The assumption (G) is equivalent to $E$ having good reduction at $p$. From
now on, we assume in addition that $E$ has {\sl ordinary} reduction at $p$,
which is equivalent to

$$
p\nd a_p.
\leqno (Ord)
$$
This assumption implies that the local Euler factor of $E$ at $p$
factorizes as

$$
1 - a_p X + pX^2 = (1 - \alpha X)(1 - \beta X),
$$
where $\alpha, \beta\in\zp$ satisfy $\ord_p(\alpha) = 0$, $\ord_p(\beta) = 1$.

Define, for $n\geq n_0+1$,

$$
y_n = x_n\otimes\alpha^{1-n} - x_{n-1}\otimes\alpha^{-n} \in E(\kn)\otimes\zp.
$$
Then

$$
\Tr_{K_{n+1}/\kn}(y_{n+1}) = y_n \qquad\qquad (n\geq n_0+1),
$$
i.e. $y = (y_n)$ is an element of the projective limit

$$
\invlim_{n > n_0} \(E(\kn)\otimes\zp\) = \invlim_n \(E(\kn)\otimes\zp\)
\subseteq \invlim_n X(E/\kn) =: \x.
$$
We shall also be interested in the inductive limit

$$
\s = \dirlim_n S(E/\kn).
$$
Both $\x$ and $\s$ are $\Lambda$-modules (of finite and co-finite type,
respectively).

\vskip24pt
\centerline{\bf 2. Iwasawa theory}

\vskip24pt\noindent
In this section we recall basic results of Iwasawa theory of elliptic
curves relating the $\Lambda$-modules $\x$ and $\s$ to Selmer groups over
the fields $\kn$. The assumptions (Heeg), (G) and (Ord) are in force.

\proclaim {Lemma 2.1}. {\rm (i)} \ $(\forall n\geq 0)$ the canonical map
$$
S(E/\kn) \lo (\s)^\gn
$$
has finite kernel and cokernel.\hb
{\rm (ii)} \ There is an isomorphism of $\Lambda$-modules (of finite type)
$$
\x \iso \Hom_\Lambda(\hat\s,\Lambda)
$$
(where $\hat M$ denotes the Pontryagin dual of $M$), hence\
$\rk_\Lambda(\x) = \corank_\Lambda(\s)$.\hb
{\rm (iii)} \ $(\forall n\geq 0)$ the canonical map
$$
(\x)_\gn \lo X(E/\kn)
$$
has finite kernel.\hb
{\rm (iv)} \ $\x\iso\Lambda^r$ for some $r\geq 0$.\hb
{\rm (v)} \ $E(\ki)_\tors$ is finite.

\sn
{\it Proof.\/} (i) \ [Man], Thm. 4.5 or [Maz1], Prop. 6.4 (note that [Man],
Lemma 4.6 eliminates the need for the second assumption in [Maz1], 6.1;
however, the latter is satisfied in our situation, because of (Heeg)).
(ii) \ [PR], Lemma 5, p. 417. (iii) \ This follows from (i) and (ii)
(cf. [PR], Lemma 4, p. 415). (iv) \ The R.H.S. in (ii) is reflexive,
hence free over $\Lambda$.  (v) \ [NeSc], 2.2.

\proclaim {Lemma 2.2}. \ $y\not=0$ in $\x$.

\sn
{\it Proof.\/} If $y_n = 0$ in $E(\kn)\otimes\zp$ for all $n > n_0$,
then

$$
x_n\otimes\alpha = x_{n-1}\otimes 1
$$
in $E(\kn)\otimes\qp$ for all $n > n_0$. As both $x_n\otimes\ 1$ and
$x_{n-1}\otimes 1$ lie in $E(\kn)\otimes\Q$, but $\alpha\not\in\Q$,
we must have $x_n\otimes\ 1 = 0$ in $E(\kn)\otimes\Q$ for all $n > n_0$,
which contradicts (now proven) Mazur's conjecture 1.4.

\proclaim {Lemma 2.3}. \ $\x\iso\Lambda$.

\sn
{\it Proof.\/} (Note that the implication [Mazur's conjecture
$\Lo \x \iso\Lambda$] was proved under more restrictive assumptions
by Bertolini [Be]).

\noindent
By Lemma 2.2 there is an exact sequence of $\Lambda$-modules of finite type

$$
0 \lo \Lambda y \lo \x \lo \x/\Lambda y \lo 0.
$$
For every {\sl surjective} character

$$
\chi : \Gamma \lo \mu_{p^n}
$$
we put $c(\chi) = n$ and denote by

$$
e_\chi = {1\over {p^n}} \sum_{g\in\Gamma/\gn} \chi(g)^{-1} g
\in \qp(\mu_{p^n})[\Gamma/\gn]
$$
the corresponding idempotent. Denote by $A$ the (finite) set of integers
$n\geq 1$ such that

$$
{{\omega_n} \over {\omega_{n-1}}} \in \Supp_\Lambda((\x/\Lambda y)_\tors)
$$
(where $\omega_n = \gamma^{p^n} - 1$, as usual). If $c(\chi)\not\in A$,
then

$$
e_\chi\(\qp(\mu_{p^n})\cdot y\ {\rm mod\ }\omega_n\x\) \not=0
$$
in $(\x)_\gn\otimes_\zp \qp(\mu_{p^n})$, hence also

$$
e_\chi\(\qp(\mu_{p^n})\cdot y_n\) \not=0
$$
in $X(E/\kn)\otimes_\zp \qp(\mu_{p^n})$, by Lemma 2.1(iii).

According to a mild generalization of the main result of [Be-Da],

$$
e_\chi\(\qp(\mu_{p^n})\cdot y_n\) \not=0 \Lo
e_\chi\(\qp(\mu_{p^n})\cdot y_n\) =
e_\chi\(X(E/\kn)\otimes_\zp \qp(\mu_{p^n})\).
$$
Putting everything together and again appealing to Lemma 2.1(iii), we
see that

$$
\rk_\zp ((\x)_\gn) \leq p^n + O(1),
$$
hence $r\leq 1$ in Lemma 2.1(iv). However, $r\not=0$ by Lemma 2.2.

\sn
{\bf 2.4} \ Recall that finite $\Lambda$-modules are also called
{\sl pseudo-null}. Denote by $(\Mod\Lambda)$ the category of all
$\Lambda$-modules and by $(\Mod\Lambda)/(ps-null)$ the category
obtained from $(\Mod\Lambda)$ by inverting all morphisms which have
pseudo-null both kernel and cokernel. This is again an abelian
category.

\proclaim {Lemma 2.5}. \ In $(\Mod\Lambda)/(ps-null)$ there is an exact
sequence
$$
0 \lo Y\oplus Y \oplus Z \lo \hat\s \lo \Lambda \lo 0
$$
and an isomorphism

$$
Z \iso \bigoplus_{i=1}^k (\Lambda/p^{m_i}\Lambda).
$$

\sn
{\it Proof.\/} Combining Lemma 2.3 with Lemma 2.1(ii) we get an exact
sequence in $(\Mod\Lambda)/(ps-null)$

$$
0 \lo (\hat\s)_\tors \lo \hat\s\lo \Lambda \lo 0.
$$
The duality results of [Ne] imply that

$$
(\hat\s)_\tors \iso Y\oplus Y \oplus Z
$$
in $(\Mod\Lambda)/(ps-null)$, with $Z$ as in the statement of the Lemma
(in fact, the condition (Heeg) implies that $Z$ itself is pseudo-null,
at least if $p > 2$; however, we do not need this fact).

\vskip24pt
\centerline{\bf 3. Main Results}

\vskip24pt

\proclaim Theorem B. Under the assumptions (Heeg) and $p\nd N a_p$ we have
$$
\sp K \equiv 1 \equiv \ran K \pmod 2 .
$$

\sn
{\it Proof.\/} First of all, (Heeg) implies that $\ran K$ is odd. Lemma 2.5
together with Lemma 2.1(i) give an exact sequence

$$
0 \lo (Y_\Gamma \oplus Y_\Gamma)\otimes_\zp \qp \lo \Hom_\zp(X(E/K),\qp)
\lo \qp \lo 0,
$$
hence

$$
\sp K = 1 + 2\, \rk_\zp (Y_\Gamma) \equiv 1 \pmod 2 .
$$

\proclaim Theorem A. Let $E$ be an elliptic curve over $\Q$ with good
ordinary reduction at $p$. Then
$$
\ran \Q \equiv \sp \Q \pmod 2 .
$$

\sn
{\it Proof.\/} For every $K = \Q(\sqrt D)$ satisfying (Heeg) we have

$$
\eqalign{
\sp K &= \sp \Q + s_p(E_D/\Q)\cr
\ran K &= \ran \Q + r_{an}(E_D/\Q),\cr}
$$
where

$$
E_D : Dy^2 = x^3 + Ax + B
$$
is the quadratic twist of $E$ over $K$. We distinguish two cases:

\sn
(I) \ $\ran \Q$ is odd.

\sn
According to [Wa] there exists $K$ satisfying (Heeg)
such that

$$
r_{an}(E_D/\Q) = 0.
$$
Results of Kolyvagin [Ko] then imply $s_p(E_D/\Q) = 0$, hence

$$
\sp \Q = \sp K \equiv \ran K = \ran \Q \pmod 2
$$
by Theorem B.

\sn 
(II) \ $\ran \Q$ is even.

\sn
Choose any $K$ satisfying (Heeg) and $p\nd D$. Applying
the result of Case (I) to $E_D$ (which has good ordinary reduction at $p$),
we obtain

$$
s_p(E_D/\Q) \equiv r_{an}(E_D/\Q) \pmod 2 ,
$$
hence

$$
\sp \Q \equiv 1 - s_p(E_D/\Q) \equiv 1 - r_{an}(E_D/\Q) \equiv \ran \Q
\pmod 2 ,
$$
again using Theorem B.

\vskip24pt
\centerline{\bf 4. Higher dimensional quotients of $J_0(N)$}

\vskip24pt\noindent
Results of Section 3 can be generalized as follows.

\sn
{\bf 4.1} \ Let $f = \sum_{n\geq 1} a_n q^n \in S_2(\Gamma_0(N))$ be
a normalized ($a_1 = 1$) newform. Put $F(f) = \Q(a_1, a_2, \cdots)$;
this is a totally real number field. Let $A$ be a quotient abelian variety
of $J_0(N)$ corresponding to $f$; it has dimension $[F(f):\Q]$, is defined
over $\Q$ and is unique up to isogeny. One has an embedding
$\iota : F(f) \ho \End(A)\otimes\Q$; for each $n\geq 1$, $\iota(a_n)$ is
induced by the Hecke correspondence $T(n) : J_0(N) \lo J_0(N)$ (with respect
to the Albanese functoriality).

\sn
{\bf 4.2} \ For each prime number $p$, $V_p(A) = T_p(A)\otimes_\zp\qp$
is a free module of rank $2$ over $F(f)\otimes_\Q\qp = \bigoplus_{\wp\vert p}
F(f)_\wp$. Fix a prime $\wp$ above $p$ in $F(f)$; then the
$F(f)_\wp$-component $V_\wp(A)$ of $V_p(A)$ defines a two-dimensional
Galois representation

$$
G_\Q = G(\qb/\Q) \lo \Aut_{F(f)_\wp}(V_\wp(A)) \iso GL_2(F(f)_\wp)
$$
which is unramified outside $Np$ and satisfies

$$
\det(1 - X\cdot Fr_{\rm geom}(\ell) \vert V_\wp(A)) =
1 - a_\ell \ell^{-1} X + \ell^{-1} X^2
$$
for all prime numbers $\ell\nd Np$.

\sn
{\bf 4.3} \ Let $\O_\wp$ be the ring of integers of $F(f)_\wp$. For every
number field $F$, denote by $S_\wp(A/F)$ the $\wp$-primary part of the Selmer
group $S_p(A/F) = \dirlim_n S(A/F, p^n)$. This is an $\O_\wp$-module of
cofinite type; denote by $s_\wp(A/F)$ its $\O_\wp$-corank. Note that
$S_\wp(A/F)$ depends only on the $G(\qb/F)$-module $V_\wp(A)/T_\wp(A)$,
where $T_\wp(A)$ is the $\O_\wp$-component of $T_p(A)$, as it coincides
with a Bloch-Kato Selmer group

$$
S_\wp(A/F) = H^1_f(F,V_\wp(A)/T_\wp(A)).
$$
The role of $\ran F$ is played by the order of vanishing

$$
r_{an}(f,F) = \ord_{s=1} L(f\otimes F,s).
$$
Under the assumptions (Heeg), (G) and

$$
\ord_\wp(a_p) = 0,
\leqno (Ord^\prime)
$$
the arguments in Section 2 go through for the $\O_\wp$-modules $S_\wp(A/K_n)$
(Mazur's control Theorem has to be replaced by a purely cohomological
``control theorem"; see [Gre]). Similarly, all of the arguments of
Section 3 work, if we replace the reference [Ko] by [KoLo]. The final
results are the following.

\proclaim Theorem B'. Under the assumptions (Heeg), $p\nd N$ and
$\ord_\wp(a_p) = 0$ we have
$$
s_\wp(A/K) \equiv 1 \equiv r_{an}(f,K) \pmod 2 .
$$

\proclaim Theorem A'. Assume that $p\nd N$ and $\ord_\wp(a_p) = 0$. Then
$$
\ord_{s=1} L(f,s) \equiv s_\wp(A/\Q) \pmod 2 .
$$

\sn
More general results can be deduced by applying the techniques of [NePl].
This will be discussed in a separate publication.

\vskip16pt
\centerline{\bf References}
\bn

\item{[Be]}
M. Bertolini, {\it Selmer groups and Heegner points in anticyclotomic
$\Z_p$-extensions}, Compositio Math. {\bf 99} (1995), no. 2, 153--182. 

\item{[BeDa]}
M. Bertolini, H. Darmon, {\it Kolyvagin's descent and Mordell-Weil groups
over ring class fields}, J. Reine Angew. Math. {\bf 412} (1990), 63--74.

\item{[Bi]}
B. Birch, {\it Heegner points of elliptic curves}, In: Symposia Mathematica,
Vol. XV (Convegno di Strutture in Corpi Algebrici, INDAM, Rome, 1973),
pp. 441--445. Academic Press, London, 1975.

\item{[BCDT]}
C. Breuil, B. Conrad, F. Diamond, R. Taylor, {\it On the modularity of
elliptic curves over $\Q$}, preprint.

\item{[Co]}
C. Cornut, {\it Mazur's conjecture on higher Heegner points}, preprint.

\item{[Gre]}
R. Greenberg, {\it Iwasawa theory for $p$-adic representations}, In:
Algebraic number theory, 97--137, Adv. Stud. Pure Math. {\bf 17}, 
Academic Press, Boston, MA, 1989.

\item{[Gro]}
B.H. Gross, {\it Heegner points on $X_0(N)$}, In: Modular forms (Durham,
1983), 87--105, Ellis Horwood Ser. Math. Appl., Horwood, Chichester, 1984.

\item{[Ko]}
V.A. Kolyvagin, {\it Euler systems}, In: The Grothendieck Festschrift,
Vol. II, 435--483, Progress in Math. {\bf 87}, Birkh\"auser Boston,
Boston, MA, 1990.

\item{[KoLo]}
V.A. Kolyvagin, D.Yu. Logachev, {\it Finiteness of the Shafarevich-Tate
group and the group of rational points for some modular abelian varieties}
(Russian), Algebra i Analiz {\bf 1} (1989), No. 5, 171--196. English
translation: Leningrad Math. J. {\bf 1} (1990), No. 5, 1229--1253.

\item{[Man]}
Yu.I. Manin, {\it Cyclotomic fields and modular curves}, Uspehi Mat. Nauk
{\bf 26} (1971), no. 6(162), 7--71. (Russian) 

\item{[Maz1]}
B. Mazur, {\it Rational points of abelian varieties with values in towers
of number fields}, Invent. Math. {\bf 18} (1972), 183--266.

\item{[Maz2]}
B. Mazur, {\it Modular curves and arithmetic}, In: Proceedings of the
ICM 1983 (Warsaw), Vol. 1, 185--211, PWN, Warsaw, 1984. 

\item{[Ne]}
J. Nekov\'a\v r, {\it Selmer complexes}, in preparation.

\item{[NePl]}
J. Nekov\'a\v r, A. Plater, {\it On the parity of ranks of Selmer groups},
Asian J. Math. {\bf 4} (2000), No. 2, 437--498.

\item{[NeSc]}
J. Nekov\'a\v r, N. Schappacher, {\it On the asymptotic behaviour of
Heegner points}, Turkish J. of Math. {\bf 23} (1999), No. 4, 549--556.

\item{[PR]}
B. Perrin-Riou, {\it Fonctions $L$ $p$-adiques, th\'eorie d'Iwasawa et
points de Heegner}, Bull. Soc. Math. France {\bf 115} (1987), No. 4, 399--456.

\item{[Va1]}
V. Vatsal, {\it Uniform distribution of Heegner points}, preprint.

\item{[Va2]}
V. Vatsal, {\it Special values of anticyclotomic $L$-functions}, preprint.

\item{[Wa]}
J.-L. Waldspurger, {\it Correspondances de Shimura}, In: Proceedings of the
ICM 1983 (Warsaw), Vol. 1, 525--531, PWN, Warsaw, 1984.

\vskip16pt
\noindent
Department of Pure Mathematics and Mathematical Statistics

\noindent
University of Cambridge

\noindent
Centre for Mathematical Sciences

\noindent
Wilberforce Road

\noindent
Cambridge CB3 0WB

\noindent
UK

\noindent
{\tt nekovar@dpmms.cam.ac.uk}

\bye